\documentclass[12pt]{amsart}
\usepackage{bm}

\newcommand{\newsection}[1]{\setcounter{equation}{0} \section{#1}}
\newcommand{\bea}{\begin{eqnarray}}
\newcommand{\eea}{\end{eqnarray}}

\newcommand{\cle}{\mathcal{E}}
\newcommand{\clf}{\mathcal{F}}

\newcommand{\clh}{\mathcal{H}}
\newcommand{\clk}{\mathcal{K}}
\newcommand{\cll}{\mathcal{L}}
\newcommand{\clm}{\mathcal{M}}
\newcommand{\cln}{\mathcal{N}}
\newcommand{\clo}{\mathcal{O}}

\newcommand{\clr}{\mathcal{R}}
\newcommand{\cls}{\mathcal{S}}

\newcommand{\raro}{\rightarrow}

\def \qed {\hfill \vrule height6pt width 6pt depth 0pt}
\def\textmatrix#1&#2\\#3&#4\\{\bigl({#1 \atop #3}\ {#2 \atop #4}\bigr)}
\def\dispmatrix#1&#2\\#3&#4\\{\left({#1 \atop #3}\ {#2 \atop #4}\right)}
\newcommand{\be}{\begin{equation}}
\newcommand{\ee}{\end{equation}}
\newcommand{\ben}{\begin{eqnarray*}}
\newcommand{\een}{\end{eqnarray*}}

\newcommand{\NI}{\noindent}

\newcommand{\bi}{\begin{itemize}}
\newcommand{\ei}{\end{itemize}}

\newcommand\la{{\langle }}
\newcommand\ra{{\rangle}}

\newtheorem{Theorem}{\sc Theorem}[section]
\newtheorem{Lemma}[Theorem]{\sc Lemma}
\newtheorem{Proposition}[Theorem]{\sc Proposition}
\newtheorem{Corollary}[Theorem]{\sc Corollary}
\newtheorem{Definition}[Theorem]{\sc Definition}
\newtheorem{Example}[Theorem]{\sc Example}
\newtheorem{Remark}[Theorem]{\sc Remark}
\newtheorem{Note}[Theorem]{\sc Note}
\newtheorem{Question}{\sc Question}
\newtheorem{Statement}{\sc Statement}
\newtheorem{ass}[Theorem]{\sc Assumption}
\newcommand{\bt}{\begin{Theorem}}
\def\beginlem{\begin{Lemma}}
\def\beginprop{\begin{Proposition}}
\def\begincor{\begin{Corollary}}
\def\begindef{\begin{Definition}}
\def\beginexamp{\begin{Example}}
\def\beginrem{\begin{Remark}}
\def\beginq{\begin{Question}}
\def\beginass{\begin{ass}}
\def\beginnote{\begin{Note}}
\newcommand{\et}{\end{Theorem}}
\def\endlem{\end{Lemma}}
\def\endprop{\end{Proposition}}
\def\endcor{\end{Corollary}}
\def\enddef{\end{Definition}}
\def\endexamp{\end{Example}}
\def\endrem{\end{Remark}}
\def\endq{\end{Question}}
\def\endass{\end{ass}}
\def\endnote{\end{Note}}

\begin{document}

\title[Resolutions of Hilbert Modules and Similarity]{Resolutions of Hilbert Modules and Similarity}
\author[Douglas]{Ronald G. Douglas}
\author[Foias]{Ciprian Foias}
\author[Sarkar]{Jaydeb Sarkar}

\address{Texas A \& M University, College Station, Texas 77843, USA}
\email{rdouglas@math.tamu.edu}
\email{jsarkar@math.tamu.edu}

\address[Present Address of Jaydeb Sarkar]{Department of Mathematics, The University of Texas at San Antonio, San Antonio, TX 78249, USA}
\email{jaydeb.sarkar@utsa.edu}

\thanks{This research was partially supported by a grant from the National Science Foundation.}



\begin{abstract}Let $H^2_m$ be the Drury-Arveson (DA) module which is the reproducing kernel Hilbert space with the kernel function $(z, w) \in \mathbb{B}^m \times \mathbb{B}^m \raro (1 - \sum\limits_{i=1}^{m} z_i \bar{w}_i)^{-1}$. We investigate for which multipliers  $\theta : \mathbb{B}^m \raro \cll(\cle, \cle_*)$ with $\mbox{ran}\, M_{\theta}$ closed, the quotient module $\clh_{\theta}$, given by $$ \cdots \longrightarrow H^2_m \otimes \cle \stackrel{M_{\theta}}{\longrightarrow} H^2_m \otimes \cle_* \stackrel{\pi_{\theta}} \longrightarrow \clh_{\theta} \longrightarrow 0,$$ is similar to $H^2_m \otimes \clf$ for some Hilbert space $\clf$. Here $M_{\theta}$ is the corresponding multiplication operator in $\cll(H^2_m \otimes \cle, H^2_m \otimes \cle_*)$ for Hilbert spaces $\cle$ and $\cle_*$ and $\clh_{\theta}$ is the quotient module $(H^2_m \otimes \cle_*)/ M_{\theta}(H^2_m \otimes \cle)$, and $\pi_{\theta}$ is the quotient map. We show that a necessary condition is the existence of a multiplier $\psi$ in $\clm(\cle_*, \cle)$ such that $$\theta \psi \theta = \theta.$$ Moreover, we show that the converse is equivalent to a structure theorem for complemented submodules of $H^2_m \otimes \cle$ for a Hilbert space $\cle$, which is valid for the case of $m=1$.
The latter result generalizes a known theorem on similarity to the unilateral shift, but the above statement is new. Further, we show that a {\it finite} resolution of DA-modules of arbitrary multiplicity using partially isometric module maps must be trivial. Finally, we discuss the analogous questions when the underlying operator $m$-tuple (or algebra) is not necessarily commuting (or commutative). In this case the converse to the similarity result is always valid.
\end{abstract}

\maketitle

\newsection{Introduction}
A well known result in operator theory (see \cite{NF63} and \cite{NF73}) states that the contraction operator given by a canonical model is similar to a unilateral shift of some multiplicity if and only if its characteristic function has a left inverse. Various approaches to this one-variable result have been given (cf. \cite{Ni}) but a new one is given in this paper which uses the commutant lifting theorem (CLT). In particular, the proof does not involve, at least explicitly, the geometry of the dilation space for the contraction.

The Drury-Arveson (DA) space $H^2_m$ (see  \cite{D}, \cite{MV}, \cite{A}) has been intensively studied by many researchers over the past few decades. In particular, the CLT has been extended to this space with a few necessary changes. Using the CLT, we extend to the DA space one direction of the one variable result on the similarity of quotient modules of the Hardy space on the unit disk. We show that the converse is equivalent to the assertion that each complemented submodule of $H^2_m \otimes \cle$ for a Hilbert space $\cle$ is isomorphic to $H^2_m \otimes \cle_*$ for some Hilbert space $\cle_*$. Of course this result follows trivially from the Beurling-Lax-Halmos theorem (BLHT) in case $m=1$. (Actually, for $m=1$ the submodule is isometrically isomorphic to $H^2_1 \otimes \cle_* = H^2(\mathbb{D}) \otimes \cle_*$.)

In Section 2, we recall some definitions and results in multivariable operator theory. In the next section, we consider a characterization of those pure co-spherically contractive Hilbert modules similar to the DA-module of some multiplicity. Using the representation of submodules of the DA-module by inner multipliers \cite{MT}, we are able to obtain the characterization in terms of inner multiplier associated with a given quotient Hilbert module and the regular inverse of that multiplier.

The quotient modules described above are the simplest case of a resolution by DA-modules for which the connecting maps are all partial isometries or inner multipliers (see \cite{G}). More precisely, using the results of Arveson \cite{A} and Muller and Vasilescu \cite{MV} and McCullough and Trent \cite{MT},  for a given pure co-spherical contractive Hilbert module one can obtain an inner resolution. In \cite{A04}, Arveson suggested that the inner resolution might not terminate as resolutions do in the algebraic context. In this paper we show that the only isometric inner multiplier, $V : H^2_m \otimes \cle \raro H^2_m \otimes \cle_*$ for Hilbert spaces $\cle$ and $\cle_*$, is the trivial one determined by an isometric operator $V_0 : 1 \otimes \cle \raro 1 \otimes \cle_*$. As a consequence, we show that all finite inner resolutions are trivial in a sense that will be explained in Section 4.

In Section 5, we are able to apply essentially the same proofs to the non-commutative case to obtain an analogous result, except here we need the noncommutative analogue of the BLHT due to Popescu (\cite{P06}, \cite{P89}). More precisely, we show that a quotient of the Fock Hilbert space, $F^2_m \otimes \cle$, for some Hilbert space $\cle$, by the range of a multi-analytic map $\Theta$ is similar to $F^2_m \otimes \clf$ for some Hilbert space $\clf$ if and only if $\Theta$ has a multi-analytic regular inverse.

In a concluding section we indicate that many of these results can be extended to complete Nevanlinna-Pick kernel Hilbert spaces and to other Hilbert modules for which the CLT holds.

\vspace{0.1in}

\NI\thanks{\textsf{Acknowledgement:} The authors wish to thank the referee for a careful reading of the
manuscript and useful remarks which led to an improved paper.}

\vspace{0.3in}

\newsection{Preliminaries}
We consider two cases, the first one in which the operators commute, or for which the algebra is $\mathbb{C}[z_1, \ldots, z_m]$ and hence commutative, and the second in which the operators are not assumed to commute or the algebra is $\mathbb{F}[Z_1, \ldots, Z_m]$. We begin with the commutative case.

Let $\{T_1, \ldots, T_m\}$ be a commuting $m$-tuple of bounded linear operators on a Hilbert space $\clh$; that is, $[T_i, T_j] = T_i T_j - T_j T_i = 0$ for $i, j = 1, \ldots, m$. A Hilbert module $\clh$ over the polynomial algebra $\mathbb{C}[z_1, \ldots, z_m]$ of $m$ commuting variables is defined so that the module multiplication $\mathbb{C}[z_1, \ldots, z_m] \times \clh \raro \clh$ is defined by $$p(z_1, \ldots, z_m) \cdot h = p(T_1, \ldots, T_m) h,$$ where $p(z_1, \ldots, z_m) \in \mathbb{C}[z_1, \ldots, z_m]$ and $h \in \clh$. We denote  by $M_1, \ldots, M_m$ the operators defined to be module multiplication by the coordinate functions. More precisely, $$M_i h = z_i \cdot h = T_i h,     \quad \quad (h \in \clh, i = 1, \ldots, m).$$

All submodules in this paper are assumed to be closed in the norm topology.

A Hilbert module over $\mathbb{C}[z_1, \ldots, z_m]$ is said to be {\it co-spherically contractive}, or define a  {\it row contraction}, if $$\|\sum_{i=1}^{m} M_i h_i\|^2 \leq \sum_{i=1}^{m} \|h_i\|^2, \quad (h_1, \ldots, h_m \in \clh),$$
or, equivalently, if $$\sum_{i=1}^{m} M_i M_i^* \leq I_{\clh}.$$

Natural examples of co-spherically contractive Hilbert modules over $\mathbb{C}[z_1, \ldots, z_m]$ are the DA-module, the Hardy module and the Bergman module, all defined on the unit ball $\mathbb{B}^m$ in $\mathbb{C}^m$. These are all reproducing kernel Hilbert spaces over $\mathbb{B}^m$ and, among them, the DA-module plays the key role for the class of co-spherically contractive Hilbert modules over $\mathbb{C}[z_1, \ldots, z_m]$. In order to be more precise, we briefly recall that a scalar reproducing kernel $K$ on a set $X$ is a function $K : X \times X \raro \mathbb{C}$ which satisfies $$\sum_{i, j = 1}^{l} \bar{c}_i c_j K(x_i, x_j) > 0,$$ for $x_1,\ldots,x_l\in X$, $c_1,\ldots,c_l \in \mathbb{C}$ with not all $c_i$ zero and $l \in \mathbb{N}$. The reproducing kernel Hilbert space $\clh_K$, corresponding to the kernel $K$, is the Hilbert space of functions defined on $X$ with the following reproducing property $$f(x) = \langle f, K_x \rangle, \; f \in \clh_K,$$where for each $x \in X$, $K_x : X \raro \mathbb{C}$ is the vector in $\clh_K$ defined by $K_x(w) = K(w, x)$, $w \in X$. The DA-module $H^2_m$ is the reproducing kernel Hilbert space corresponding to the kernel $K : \mathbb{B}^m \times \mathbb{B}^m \raro \mathbb{C}$ defined by $$K(z, w) = (1 - \sum_{i=1}^{m} z_i \bar{w}_i)^{-1}, \; (z, w) \in \mathbb{B}^m \times \mathbb{B}^m.$$

We identify the Hilbert tensor product $H^2_m \otimes \cle$ with the $\cle$-valued $H^2_m$ space $H^2_m(\cle)$ or the $\cll(\cle)$-valued reproducing kernel Hilbert space with the kernel $(z, w) \mapsto (1 - \sum\limits_{i=1}^{m} z_i \bar{w}_i)^{-1} I_{\cle}$. Consequently,
$$H^2_m \otimes \cle = \{ f \in \clo (\mathbb{B}^m, \cle) : f(z) = \sum_{\bm{k} \in \mathbb{N}^m} a_{\bm{k}} z^{\bm{k}}, a_{\bm{k}} \in \cle, \|f\|^2 : = \sum_{\bm{k} \in \mathbb{N}^m} \frac{ \| a_{\bm{k}} \|^2}{\gamma_{\bm{k}}} < \infty \},$$
where $\clo(\mathbb{B}^m, \cle)$ is the space of $\cle$-valued holomorphic functions on $\mathbb{B}^m$, $\bm{k} = (k_1, \ldots, k_m)$ and $\gamma_{\bm{k}} = \frac{(k_1 + \cdots + k_m)!}{k_1 ! \cdots k_m!}$ are the multinomial coefficients.
A function $\varphi \in \clo(\mathbb{B}^m, \cll(\cle, \cle_*))$ is said to be a \textit{multiplier} if $\varphi f \in H^2_m \otimes \cle_* = H^2_m(\cle_*)$ for all $f \in H^2_m \otimes \cle = H^2_m(\cle)$. By the closed graph theorem, such a multiplier $\varphi$ defines a bounded module map $$M_{\varphi} : H^2_m \otimes \cle \raro H^2_m \otimes \cle_*, \quad M_{\varphi} f = \varphi f, \,f \in H^2_m \otimes \cle.$$Equivalently, we can consider $\varphi \in \clo (\mathbb{B}^m, \cll(\cle, \cle_*))$ for which $M_{\varphi}$ defines a bounded operator from $H^2_m \otimes \cle$ to $H^2_m \otimes \cle_*$. The set of all such bounded multipliers $\varphi \in \clo (\mathbb{B}^m, \cll(\cle, \cle_*))$ will be denoted by $\clm(\cle, \cle_*)$. A multiplier $\varphi \in \clm(\cle, \cle_*)$ is said to be \textit{inner} if $M_{\varphi}$ is a partial isometry in $\cll(H^2_m \otimes \cle, H^2_m \otimes \cle_*)$.

We recall an analogue of the CLT due to Ball-Trent-Vinnikov (Theorem 5.1 in \cite{BTV}) on DA-modules which will be used to prove some of the main results of this paper.

\begin{Theorem}\label{CLT} \textsf{(Ball-Trent-Vinnikov)}
Let $\cln$ and $\cln_*$ be quotient modules of $H^2_m \otimes \cle$ and $H^2_m \otimes \cle_*$ for some Hilbert spaces $\cle$ and $\cle_*$, respectively. If $X : \cln \raro \cln_*$ is a bounded module map, that is,   $$X P_{\cln} (M_{z_i} \otimes I_{\cle})|_{\cln} = P_{\cln_*} (M_{z_i} \otimes I_{\cle_*})|_{\cln_*} X,$$ for $i = 1, \ldots, m$, then there exists a multiplier $\varphi \in \clm(\cle, \cle_*)$ such that\\
\indent (i) $\|X\| = \|M_{\varphi}\|$ and

(ii) $P_{\cln_*} M_{\varphi} = X.$

\NI In the language of Hilbert modules, one has $\pi_{\cln_*} M_{\varphi} = X \pi_{\cln}$, where $\pi_{\cln}$ and $\pi_{\cln_*}$ are the quotient maps.
\end{Theorem}

The above statement of the CLT for $\mathbb{C}[z_1,\ldots, z_m]$ is due to Ball-Trent-Vinnikov as indicated. However, Popescu pointed out that the result follows from its noncommutative analogue established earlier by him in \cite{P89T,P89}. A more recent paper on this topic is due to Davidson and Le (\cite{DL}).

We now recall the notion of pureness for a co-spherically contractive Hilbert module $\clh$ over $\mathbb{C}[z_1, \ldots, z_m]$. Define the completely positive map $$P_{\clh} : \cll(\clh) \rightarrow \cll(\clh)$$ by $$P_{\clh} (A) = \sum_{i=1}^{m} M_i A M^*_i, \quad \quad \quad \quad (A \in \cll(\clh)).$$ Now $$I_{\clh} \geq P_{\clh} (I_{\clh}) \geq P^2_{\clh} (I_{\clh}) \geq \cdots \geq P^l_{\clh} (I_{\clh}) \geq \cdots \geq 0,$$ so that $$P_{\infty} = \mbox{SOT} - \mbox{lim}_{l \raro \infty} P_{\clh}^l (I_{\clh})$$ exists and $0 \leq P_{\infty} \leq I_{\clh}$. The Hilbert module $\clh$ is said to be \textit{pure} if $$P_{\infty} = 0.$$

A canonical example of a pure co-spherically contractive Hilbert module over $\mathbb{C}[z_1, \ldots, z_m]$ is the DA-module $H^2_m \otimes \clf$, where $\clf$ is a Hilbert space. Moreover, quotients of DA-modules characterize all pure co-spherically contractive Hilbert modules.

\begin{Theorem}\label{AMVMT}
Let $\clh$ be a pure co-spherically contractive Hilbert module. Then

(i)\title{(\textsf{Arveson} \cite{A}, \textsf{Muller-Vasilescu} \cite{MV})} $\clh$ is isometrically isomorphic with a quotient module $(H^2_m \otimes \cle_*)\,/ \,\cls$, where $\cle$ is a Hilbert space and $\cls$ is a submodule of $H^2_m \otimes \cle_*$.

(ii) \title{(\textsf{McCullough-Trent} \cite{MT})} If $\cls$ is a submodule of $H^2_m \otimes \cle_*$, then there exists a multiplier $\theta \in \clm(\cle, \cle_*)$ for some Hilbert space $\cle$ such that $M_{\theta}$ is inner and $\cls = M_{\theta} (H^2_m \otimes \cle)$.

\NI Therefore, $\clh$ is isometrically isomorphic to $(H^2_m \otimes \cle_*)\,/\,\mbox{ran} \,M_{\theta}$ for some inner multiplier $\theta \in \clm(\cle, \cle_*)$ and Hilbert space $\cle$.
\end{Theorem}

Note that in the statement of Theorem \ref{CLT}, Ball-Trent-Vinnikov \cite{BTV} made the additional assumption that the submodules $\cln^{\perp}$ and $\cln_*^{\perp}$ are invariant under the scalar multipliers. However, that this condition is redundant follows from part (ii) of Theorem \ref{AMVMT} above due to McCullough-Trent.

We now consider some preliminaries for the case of noncommuting operators. Let $\mathbb{F}^+_m$ denote the free semigroup with the $m$ generators $g_1, \ldots, g_m$ and let $F^2_m$ be the full Fock space of $m$ variables, which is a Hilbert space. More precisely, if we let $\{e_1, \ldots, e_m\}$ be the standard orthonormal basis of $\mathbb{C}^m$, then $$F^2_m = \bigoplus_{k \geq 0} (\mathbb{C}^m)^{\otimes k},$$ where $(\mathbb{C}^m)^{\otimes 0} = \mathbb{C}.$ The creation, or left shift, operators $S_1, \ldots, S_m$ on $F^2_m$ are defined by $$S_i f = e_i \otimes f,$$ for all $f$ in $F^2_m$ and $i = 1, \ldots, m$.

Let $\{T_1, \ldots, T_m\}$ be $m$ bounded linear operators on a Hilbert space $\clk$ which are not necessarily commuting. One can make $\clk$ into a Hilbert module over the algebra of polynomials $\mathbb{F}[Z_1, \ldots Z_m]$, in $m$ noncommuting variables, as follows: $$\mathbb{F}[Z_1, \ldots Z_m] \times \clk \raro \clk, \quad p(Z_1, \ldots, Z_m) \cdot h \mapsto p(T_1, \ldots, T_m) h, \, h \in \clk.$$ The module $\clk$ over $\mathbb{F}[Z_1, \ldots, Z_m]$ is said to be co-spherically contractive if the row operator given by module multiplication by the coordinate functions is a contraction.

A bounded linear operator $\Theta \in \cll(F^2_m \otimes \cle, F^2_m \otimes \cle_*)$, for some Hilbert spaces $\cle$ and $\cle_*$, is said to be a multi-analytic operator if it is a module map; that is, if  $$\Theta (S_i \otimes I_{\cle}) = (S_i \otimes I_{\cle_*}) \Theta, \; \; i = 1, \ldots, m.$$ Given a multi-analytic operator $\Theta$ as above, one can define a bounded linear operator $\theta : \cle \rightarrow F^2_m \otimes \cle_*$ by $$\theta x = \Theta (1 \otimes x) \quad \quad \quad \quad (x \in \cle).$$
In this correspondence of $\Theta$ and $\theta$, each uniquely determines the other. Moreover, one defines the operator coefficients $\theta_{\bm{\alpha}} \in \cll(\cle, \cle_*)$ of $\Theta$ by $$\la \theta_{\bm{\alpha}^t} x, y \ra = \langle\theta x, e_{\bm{\alpha}} \otimes y \rangle = \la \Theta (1 \otimes x), e_{\bm{\alpha}} \otimes y \ra \quad \quad \quad (x \in \cle, y \in \cle_*)$$for each $\bm{\alpha} \in \mathbb{F}^+_m$, where $\bm{\alpha}^t = g_{i_p} \cdots g_{i_1}$ for $\bm{\alpha} = g_{i_1} \cdots g_{i_p}$. It was proved by Popescu (cf. \cite{P06}) that $$\Theta = \mbox{SOT} - \mbox{lim}_{r \raro 1^-} \sum_{l=0}^{\infty} \sum_{|\bm{\alpha}| = l} r^{|\bm{\alpha}|} R^{\bm{\alpha}} \otimes \theta_{\bm{\alpha}},$$ where $R_i = U^* S_i U$ for $i=1, \ldots, m,$ are the right creation operators on $F^2_m$, $R^{\bm{\alpha}} = R_{g_{i_1}} \cdots R_{g_{i_p}}$ for $\bm{\alpha} = g_{i_1} \cdots g_{i_p}$, and $U$ is the unitary operator on $F^2_m$ defined by $U e_{\bm{\alpha}} = e_{\bm{\alpha}^t}$ for $\bm{\alpha} \in \mathbb{F}^+_m$. The set of all multi-analytic operators in $\cll(F^2_m \otimes \cle, F^2_m \otimes \cle_*)$ coincides with $R^{\infty}_m \bar{\otimes} \cll(\cle, \cle_*)$, the WOT closed algebra generated by the spatial tensor product of $R^{\infty}_m$ and $\cll(\cle, \cle_*)$, where $R^{\infty}_m = U^* F^{\infty}_m U$ and $F^{\infty}_m$ is the WOT closed algebra generated by the left creation operators, $S_1, \ldots, S_m,$ and the identity operator on $F^2_m$.

Notice that the definition of a pure co-spherically contractive Hilbert module can be extended to the noncommutative case; that is, with appropriate change of notation, the concept of a pure co-spherically contractive Hilbert module $\clk$ over $\mathbb{F}[Z_1, \ldots, Z_m]$ can be defined in a similar way. Popescu proved that any pure co-spherically contractive Hilbert module over $\mathbb{F}[Z_1, \ldots, Z_m]$ can be realized as a quotient module of $F^2_m \otimes \cle$ for some Hilbert space $\cle$ (see \cite{P89} and Theorem 2.10 and references in \cite{P06}).

\begin{Theorem}\label{Pop1} \textsf{(Popescu)}
 Given a pure co-spherically contractive Hilbert module $\clk$ over $\mathbb{F}[Z_1, \ldots, Z_m]$, there is a multi-analytic operator $\Theta$ in $\cll(F^2_m \otimes \cle, F^2_m \otimes \cle_*)$ for some Hilbert spaces $\cle$ and $\cle_*$, which is isometric such that $\clk$ is isometrically isomorphic to the quotient of $F^2_m \otimes \cle_*$ by the range of $\Theta$. Moreover, the characteristic operator function $\Theta$ is a complete unitary invariant for $\clk$.
\end{Theorem}

Finally we need to extend the analogue of the BLHT to the noncommuting setting which is due to Popescu (\cite{P89}, \cite{P06}).

\begin{Theorem}\label{BLH} \textsf{(Popescu)}
If $\cls$ is a closed subspace of $F^2_m \otimes \clf$ for some Hilbert space $\clf$, then the following are equivalent:

(i) $\cls$ is a submodule of $F^2_m \otimes \clf$.

(ii) There exists a Hilbert space $\cle$ and an (isometric) inner multi-analytic operator $\Phi : F^2_m \otimes \cle \raro F^2_m \otimes \clf$ such that $$\cls = \Phi (F^2_m \otimes \cle).$$
\end{Theorem}

\newsection{Hilbert modules over $\mathbb{C}[z_1, \ldots, z_m]$}

Let $\theta \in \clm(\cle, \cle_*)$ be a multiplier for Hilbert spaces $\cle$ and $\cle_*$ such that $M_{\theta}$ has closed range and let $\clh_{\theta}$ be the quotient module defined by the sequence $$ \cdots \longrightarrow H^2_m \otimes \cle \stackrel{M_{\theta}}{\longrightarrow} H^2_m \otimes \cle_* \stackrel{\pi_{\theta}} \longrightarrow \clh_{\theta} \longrightarrow 0,$$where $\pi_{\theta}$ is the quotient map of $H^2_m \otimes \cle_*$ onto the quotient of $H^2_m \otimes \cle_*$ by the range of $M_{\theta}$. There are several possible relationships between these objects:

\begin{Statement}\label{Q1}
The sequence splits or $\pi_{\theta}$ is right invertible; that is, there exists a module map $\sigma_{\theta} : \clh_{\theta} \raro H^2_m \otimes \cle_*$ such that $$\pi_{\theta} \sigma_{\theta} = I_{\clh_{\theta}}.$$
\end{Statement}

\begin{Statement}\label{Q2}
The multiplication operator $M_{\theta}$ has a left inverse. Equivalently, there exists a multiplier $\psi \in \clm(\cle_*, \cle)$ satisfying $\psi(z) \theta(z) = I_{\cle}$ for $z \in \mathbb{B}^m$.
\end{Statement}

\begin{Statement}\label{Q3}
The multiplier $\theta$ has a regular inverse. Equivalently, there exists a multiplier $\psi \in \clm(\cle_*, \cle)$ satisfying $$\theta(z) \psi(z) \theta(z) = \theta(z),$$for $z \in \mathbb{B}^m$.
\end{Statement}

Note that in case $\mbox{ker~}M_{\theta} = \{0\}$, Statements \ref{Q2} and \ref{Q3} are equivalent.

\begin{Statement}\label{Q3a}
The range of $M_{\theta}$ is complemented in $H^2_m \otimes \cle_*$ or there exists a submodule $\cls$ of $H^2_m \otimes \cle_*$ such that $$\mbox{ran~} M_{\theta} \stackrel{\bm{.}}{+} \cls = H^2_m \otimes \cle_*.$$
\end{Statement}

\begin{Statement}\label{Q4}
The quotient Hilbert module $\clh_{\theta}$ is similar to $H^2_m \otimes \clf$ for some Hilbert space $\clf$.
\end{Statement}

\begin{Statement}\label{Q5}
Suppose $H^2_m \otimes \cle_*$ is a skew direct sum $\cls_1 \stackrel{\bm{.}}{+} \cls_2$, where $\cle_*$ is a Hilbert space and $\cls_1$ and $\cls_2$ are submodules such that $\cls_1$ is isomorphic to $H^2_m \otimes \cle$ for some Hilbert space $\cle$. Then $\cls_2$ is isomorphic to $H^2_m \otimes \clf$ for some Hilbert space $\clf$.
\end{Statement}

Note that Statements \ref{Q3a} and \ref{Q5} imply Statement \ref{Q4} and would be the converse to Corollary \ref{TH1}.

\begin{Statement}\label{Q6}
If $\cls$ is a complemented submodule of $H^2_m \otimes \cle_*$ for some Hilbert space $\cle_*$, then $\cls$ is isomorphic to $H^2_m \otimes \clf$ for some Hilbert space $\clf$?
\end{Statement}

One can reformulate Statement \ref{Q6} in the following equivalent form.

\begin{Statement}\label{Q7}
Every complemented submodule $\cls$ of $H^2_m \otimes \cle$, for some Hilbert space $\cle$, is the range of $M_{\psi}$ for a multiplier $\psi \in \clm(\cle, \clf)$ with $\mbox{ker}\, M_{\psi} = \{0\}$ for some Hilbert space $\clf$.
\end{Statement}
Note that one could view an affirmation of this statement as a weak form of the BLHT for DA-modules.

Statement \ref{Q6} raises an important issue for Hilbert modules: are the complemented submodules $\cls \neq \{0\}$ of $\clr \otimes \mathbb{C}^n$ always isomorphic to $\clr \otimes \mathbb{C}^k$ for some $0 < k \leq n$. This is certainly not true for a general Hilbert module $\clr$. However, what if $\clr$ belongs to the class of ``locally-free'' Hilbert modules of multiplicity one which is the case for the DA-module $H^2_m$? For $m=1$, an affirmation follows trivially from the BLHT. A less obvious argument shows that the result holds for more general ``locally-free'' Hilbert modules over the unit disk such as the Bergman module. (Although the language is different, this result was proved by J. S. Fang, C. L. Jiang, X. Z. Guo, K. Ti and H. He. The study of the relationship between the eight statements in the one-variable case is close to the theme of the book by C. L. Jiang and F. Wang \cite{JW}, where details can be found.) Further, one can establish an affirmation to Statement \ref{Q5} if one assumes that the multiplier $\theta \in \clm(\cle, \cle_*)$ is holomorphic on a neighborhood of the closure of $\mathbb{B}^m$, at least if $\cle$ and $\cle_*$ are finite dimensional. However, what happens in general for ``locally-free'' Hilbert modules over $\mathbb{B}^m$, such as the DA-module, is not clear at this point.

Moreover, a necessary condition for $\cls_1$ and $\cls_2$ with $H^2_m \otimes \mathbb{C}^n = \cls_1 \stackrel{\bm{.}}{+} \cls_2$ to be isomorphic to $H^2_m \otimes \mathbb{C}^k$ and $H^2_m \otimes \mathbb{C}^{n-k}$, respectively, is the existence of a generating set $\{f_1, \ldots, f_n\}$ for $H^2_m \otimes \mathbb{C}^n$ with $\{f_1, \ldots, f_k\}$ in $\cls_1$ and $\{f_{k+1}, \ldots, f_n\}$ in $\cls_2$. Note one can view each $f_i \in \clo(\mathbb{B}^m, \mathbb{C}^n)$ for $i = 1, \ldots, n$. If one assumes in addition that the vectors $\{f_i\}$ are in $\clm(\mathbb{C}^n, \mathbb{C}^m)$, then that is also sufficient.

Note that given a complemented submodule $\cls$ of $H^2_m \otimes \cle$, that is, for some submodule $\tilde{\cls}$ one has $\cls \stackrel{\bm{.}}{+} \tilde{\cls} = H^2_m \otimes \cle$, there are many choices $\tilde{\tilde{\cls}}$ so that the skew direct sum $\cls \stackrel{\bm{.}}{+} \tilde{\tilde{\cls}}$ is isomorphic to $H^2_m \otimes \cle_*$ for some Hilbert space $\cle_*$. (Here we allow a different space $\cle_*$.) It is not clear, but seems unlikely that there exists a canonical choice of $\cle_*$ and $\tilde{\tilde{\cls}}$, in some sense, or what the ``simplest'' choice might be. Such ideas are related to the $K$-theory group introduced in \cite{JW}.

In commutative algebra one shows that a short exact sequence of modules $$ 0 \longrightarrow A \stackrel{\varphi_1}{\longrightarrow} B \stackrel{\varphi_2} \longrightarrow C \longrightarrow 0,$$ splits, or $\varphi_2$ has a right inverse, if and only if $\varphi_1$ has a left inverse. Using the closed graph theorem, we can extend this result to Hilbert modules. Moreover, with the CLT we can extend the result to the case where $\varphi_1$ has a kernel. We begin with the simpler result.

\begin{Theorem}\label{1}
Let $\theta \in \clm(\cle, \cle_*)$ be a multiplier for Hilbert spaces $\cle$ and $\cle_*$ such that $\mbox{ran}\,M_{\theta}$ is closed. Then $\mbox{ran}\,M_{\theta}$ is complemented in $H^2_m \otimes \cle_*$ if and only if there exists a module map $\sigma_{\theta} : \clh_{\theta} \raro H^2_m \otimes \cle_*$  which is a right inverse for $\pi_{\theta}$.
\end{Theorem}
\NI \textsf{Proof.} If $H^2_m \otimes \cle_* = \mbox{ran}\,M_{\theta} \stackrel{\cdot}+ \cls$ for a (closed) submodule $\cls$, then $Y = \pi_{\theta}|_{\cls}$ is one-to-one and onto. Hence $Y^{-1} : (H^2_m \otimes \cle_*)\,/\,\mbox{ran}\,M_{\theta} \raro \cls$ is bounded by the closed graph theorem and $\sigma_{\theta} = i ~Y^{-1}$ is a right inverse for $\pi_{\theta}$, where $i : \cls \raro H^2_m \otimes \cle_*$ is the inclusion map.

\NI Conversely, if there exists a right inverse $\sigma_{\theta} : \clh_{\theta} \raro H^2_m \otimes \cle_*$ for $\pi_{\theta}$, then $\sigma_{\theta} \pi_{\theta}$ is an idempotent on $H^2_m \otimes \cle_*$ such that $\cls = \mbox{ran}\,\sigma_{\theta} \pi_{\theta}$ is a complementary submodule for the closed submodule $\mbox{ran}\,M_{\theta}$ in $H^2_m \otimes \cle_*$. \qed

For a multiplier $\theta \in \clm(\cle, \cle_*)$ one could consider the quotient of $H^2_m \otimes \cle_*$ by the closure of $\mbox{ran~}M_{\theta}$. Examples in the case $m = 1$ show that the existence of a right inverse  for $\pi_{\theta}$ does not imply that $\mbox{ran}\, M_{\theta}$ is closed. Theorem \ref{1} shows that the Statements \ref{Q1} and \ref{Q3a} are equivalent but the preceding comment implies the necessity of the assumption that $\mbox{ran~}M_{\theta}$ is closed. However, we do have the following folklore result which helps clarify matters.

\begin{Remark}\label{2}
If $\theta \in \clm(\cle, \cle_*)$ for Hilbert spaces $\cle$ and $\cle_*$ and $\mbox{ran}\, M_{\theta}$ is complemented in $H^2_m \otimes \cle_*$, then $\mbox{ran}\, M_{\theta}$ is closed.
\end{Remark}

As one knows, by considering the $m= 1$ case, there is more than one multiplier $\theta \in \clm(\cle, \cle_*)$ for Hilbert spaces $\cle$ and $\cle_*$ with the same range and thus yielding the same quotient. While things are even more complicated for $m>1$, the following result using the CLT introduces some order.

\begin{Theorem}\label{incl}
Let $\theta \in \clm(\cle, \cle_*)$ be an inner multiplier for Hilbert spaces $\cle$ and $\cle_*$ and $\varphi \in \clm(\clf, \cle_*)$ for some Hilbert space $\clf$. Then there exists a multiplier $\psi \in \clm(\clf, \cle)$ such that $\varphi = \theta \psi$ if and only if $$\mbox{ran~} M_{\varphi} \subseteq \mbox{ran~}M_{\theta}.$$
\end{Theorem}
\NI \textsf{Proof.} If $\psi \in \clm(\clf, \cle)$ such that $\varphi = \theta \psi$, then $M_{\varphi} = M_{\theta} M_{\psi}$ and hence $$\mbox{ran~} M_{\varphi} = \mbox{ran~} M_{\theta} M_{\psi} \subseteq \mbox{ran~}M_{\theta}.$$

Suppose $\mbox{ran~} M_{\varphi} \subseteq \mbox{ran~}M_{\theta}$ and $\theta \in \clm(\cle, \cle_*)$ is an inner multiplier. This implies that $\mbox{ran}\, M_{\theta}$ is closed. Consider the module map $$\hat{M_{\theta}} : (H^2_m \otimes \cle)~/~\mbox{ker~}M_{\theta} \longrightarrow \mbox{ran~}M_{\theta}$$ defined by \[\hat{M_{\theta}} \gamma_{\theta} = M_{\theta},\]which is invertible since $\mbox{ran~} M_{\theta}$ is closed. Let $\gamma_{\theta} : H^2_m \otimes \cle \longrightarrow (H^2_m \otimes \cle)/~\mbox{ker~}M_{\theta}$ be the quotient module map. Set $\hat{X} = \hat{M_{\theta}}^{-1}$. Then $\hat{X} : \mbox{~ran~}M_{\theta} \raro (H^2_m \otimes \cle)/\mbox{~ker~}M_{\theta}$ is bounded by the closed graph theorem and so is $\hat{X} M_{\varphi} : H^2_m \otimes \clf \raro (H^2_m \otimes \cle)\,/\,\mbox{ker~} M_{\theta}$. Appealing to the CLT yields a multiplier $\psi \in \clm(\clf, \cle)$ so that $$\gamma_{\theta} M_{\psi} = \hat{X} M_{\varphi},$$ and hence \[M_{\theta} M_{\psi} = (\hat{M_{\theta}} \gamma_{\theta}) M_{\psi} = \hat{M_{\theta}} (\hat{X} M_{\varphi}) =  M_{\varphi},\] or $\varphi = \theta \psi$ which completes the proof. \qed

Note that the result holds for a multiplier $\theta \in \clm(\cle, \cle_*)$ for Hilbert spaces $\cle$ and $\cle_*$ so long as $\mbox{ran~}M_{\theta}$ is closed since that is all that the proof uses.

We now consider our principal result on multipliers and regular
inverses.

\begin{Theorem}\label{P1}
Let $\theta \in \clm(\cle, \cle_*)$ be a multiplier for Hilbert spaces $\cle$ and $\cle_*$. Then there exists $\psi \in \clm(\cle_*, \cle)$ such that $$M_{\theta} M_{\psi} M_{\theta} = M_{\theta}$$ if and only if $\mbox{ran}\, M_{\theta}$ is complemented in $H^2_m \otimes \cle_*$, or $$H^2_m \otimes \cle_* = \mbox{ran}\, M_{\theta} \stackrel{\cdot}+ \cls,$$ for some submodule $\cls$ of $H^2_m \otimes \cle_*$.
\end{Theorem}
\NI \textsf{Proof.} If $H^2_m \otimes \cle_* = \mbox{ran}\, M_{\theta} \stackrel{\cdot}+ \cls$ for some (closed) submodule $\cls$, then $\mbox{ran~}M_{\theta}$ is closed by Remark \ref{2}. Consider the module map \[\hat{M_{\theta}} : (H^2_m \otimes \cle)\,/ \,\mbox{ker}\, M_{\theta} \longrightarrow (H^2_m \otimes \cle_*)\,/\, \cls,\] defined by $$\hat{M_{\theta}} \gamma_{\theta} = \pi_{\cls} M_{\theta},$$ where $\gamma_{\theta} : H^2_m \otimes \cle \raro (H^2_m \otimes \cle)\,/\,\mbox{ker}\, M_{\theta}$ and $\pi_{\cls} : H^2_m \otimes \cle_* \raro (H^2_m \otimes \cle_*)\,/\,\cls$ are quotient maps. This map is one-to-one and onto and thus has a bounded inverse $\hat{X} = \hat{M_{\theta}}^{-1} : (H^2_m \otimes \cle_*)\,/ \,\cls \raro (H^2_m \otimes \cle)\,/\, \mbox{ker~}M_{\theta}$ by the closed graph theorem. Since $\hat{X}$ satisfies the hypotheses of the CLT, there exists $\psi \in \clm(\cle_*, \cle)$ such that $$\gamma_{\theta} M_{\psi} = \hat{X} \pi_{\cls}.$$ Further, $\hat{M_{\theta}} \gamma_{\theta} = \pi_{\cls} M_{\theta}$ yields \[\pi_{\cls} M_{\theta} M_{\psi} = \hat{M_{\theta}} \gamma_{\theta} M_{\psi} = \hat{M_{\theta}} \hat{X} \pi_{\cls} = \pi_{\cls},\] and therefore\[\pi_{\cls}(M_{\theta} M_{\psi} M_{\theta} - M_{\theta}) = 0.\] Since $\pi_{\cls}$ is one-to-one on $\mbox{ran}\, M_{\theta}$, it follows that $M_{\theta} M_{\psi} M_{\theta} = M_{\theta}$.

Now suppose there exists $\psi \in \clm(\cle_*, \cle)$ such that $M_{\theta} M_{\psi} M_{\theta} = M_{\theta}$. This implies  that $$(M_{\theta} M_{\psi})^2 = M_{\theta} M_{\psi},$$ and hence $M_{\theta} M_{\psi}$ is an idempotent. From the equality $M_{\theta} M_{\psi} M_{\theta} = M_{\theta}$ we obtain both that $\mbox{ran}\,M_{\theta} M_{\psi}$ contains $\mbox{ran}\,M_{\theta}$ and that $\mbox{ran}\, M_{\theta} M_{\psi}$ is contained in $\mbox{ran}\, M_{\theta}$. Therefore, $$\mbox{ran}\, M_{\theta} M_{\psi} = \mbox{ran}\, M_{\theta},$$ and $$\cls = \mbox{ran}\,(I - M_{\theta} M_{\psi}),$$ is a complementary submodule of $\mbox{ran~}M_{\theta}$ in $H^2_m \otimes \cle_*$. \qed

Note that Theorem \ref{P1} implies that Statements \ref{Q3} and \ref{Q3a} are equivalent.

\begin{Corollary}\label{TH1}
Assume $\theta \in \clm(\cle, \cle_*)$ for Hilbert spaces $\cle$ and $\cle_*$ such that $\mbox{ran}\, M_{\theta}$ is closed and $\clh_{\theta}$ is defined by \[H^2_m \otimes \cle \stackrel{M_{\theta}}\longrightarrow H^2_m \otimes \cle_* \longrightarrow \clh_{\theta} \longrightarrow 0.\] If $\clh_{\theta}$ is similar to $H^2_m \otimes \clf$ for some Hilbert space $\clf$, then the sequence splits.
\end{Corollary}
\NI \textsf{Proof.} First, assume that there exists an invertible module map $X : H^2_m \otimes \clf \raro \clh_{\theta}$, and let $\varphi \in \clm(\clf, \cle_*)$ be defined by the CLT so that $\pi_{\theta} M_{\varphi} = X$, where $\pi_{\theta} : H^2_m \otimes \cle_* \raro (H^2_m \otimes \cle_*)\,/\, \mbox{ran~}M_{\theta}$ is the quotient map. Since $X$ is invertible we have $$H^2_m \otimes \cle_* = \mbox{ran~} M_{\varphi} \stackrel{\bm{.}}{+}  \mbox{ran} M_{\theta}.$$ Thus $\mbox{ran~}M_{\theta}$ is complemented and hence it follows from the previous corollary that the sequence splits. \qed

Corollary \ref{TH1} shows that Statement \ref{Q5} implies Statement \ref{Q1}. If Statement \ref{Q5} is valid, then the converse to Corollary \ref{TH1} holds. Moreover, we see that Statement \ref{Q7} implies that $\clh_{\theta}$ is similar to $H^2_m \otimes \clf$ for some Hilbert space $\clf$ or that Statement \ref{Q3a} is valid. Finally, the following weaker converse to Corollary \ref{TH1} always holds.

\begin{Corollary}
Let $\theta \in \clm(\cle, \cle_*)$ for Hilbert spaces $\cle$ and $\cle_*$, and set $\clh_{\theta} = (H^2_m \otimes \cle_*)/ \mbox{~clos~}[\mbox{ran~}M_{\theta}]$. Then the following statements are equivalent:

\NI (i) there exists $\psi \in \clm(\cle_*, \cle)$ such that $\psi(x) \theta(z) = I_{\cle}$ for $z \in \mathbb{B}^m$, and

\NI (ii) $\mbox{ran~}M_{\theta}$ is closed, $\mbox{ker~}M_{\theta} = \{0\}$ and $\clh_{\theta}$ is similar to a complemented submodule $\cls$ of $H^2_m \otimes \cle_*$.
\end{Corollary}

\NI\textsf{Proof.} If (i) holds, then $\mbox{ran~}M_{\theta}$ is closed and $\mbox{ker~}M_{\theta} = \{0\}$. Further, $M_{\theta} M_{\psi}$ is an idempotent on $H^2_m \otimes \cle_*$ such that $\mbox{ran~} M_{\theta} M_{\psi} = \mbox{ran~}M_{\theta}$ and $\clh_{\theta}$ is isomorphic to $\cls = \mbox{ran~}(I - M_{\theta} M_{\psi}) \subseteq H^2_m \otimes \cle_*$ and $H^2_m \otimes \cle_* = \mbox{ran~}M_{\theta} \stackrel{\bm{.}}{+} \cls$ so $\cls$ is complemented.

\NI Now assume that (ii) holds and there exists an isomorphism $X : \clh_{\theta} \raro \cls \subseteq H^2_m \otimes \cle_*$, where $\cls$ is a complemented submodule of $H^2_m \otimes \cle_*$. Then $Y = X \pi_{\theta} : H^2_m \otimes \cle_* \raro H^2_m \otimes \cle_*$ is a module map and hence there exists a multiplier $\omega \in \clm(\cle_*, \cle_*)$ so that $Y = M_{\omega}$. Since $X$ is invertible, $\mbox{ran~}M_{\omega} = \cls$, which is complemented by assumption, and hence by Theorem \ref{P1} there exists $\psi \in \clm(\cle_*, \cle_*)$ such that $M_{\omega} = M_{\omega} M_{\psi} M_{\omega}$ or $M_{\omega}(I - M_{\psi} M_{\omega}) = 0$. Therefore, $$\mbox{ran~}(I - M_{\psi} M_{\omega}) \subseteq \mbox{ker~} M_{\omega} = \mbox{ker~}Y = \mbox{ker~}\pi_{\theta} = \mbox{ran~}M_{\theta}.$$ Applying Theorem \ref{incl}, we obtain $\varphi \in \clm(\cle_*, \cle)$ so that $$I - M_{\psi} M_{\omega} = M_{\theta} M_{\varphi}.$$ Thus using $M_{\omega} M_{\theta} = 0$ we see that \[M_{\theta} M_{\varphi} M_{\theta} = (I - M_{\psi} M_{\omega}) M_{\theta} = M_{\theta},\] or \[ M_{\theta} = M_{\theta} M_{\varphi} M_{\theta}.\] Since $\mbox{ker~}M_{\theta} = \{0\}$, we have $$M_{\varphi} M_{\theta} = I_{H^2_m \otimes \cle},$$which completes the proof. \qed

Combining Theorem \ref{P1} and Corollary \ref{TH1} yields our main result in the commutative setting for similarity.

\begin{Corollary}\label{main-cor}
Given $\theta \in \clm(\cle, \cle_*)$ for Hilbert spaces $\cle$ and $\cle_*$ such that $\mbox{ran~}M_{\theta}$ is closed, consider the quotient Hilbert module $\clh_{\theta}$ defined as above. If $\clh_{\theta}$ is similar to $H^2_m \otimes \clf$ for some Hilbert space $\clf$, then there exists a multiplier $\psi \in \clm(\cle_*, \cle)$ satisfying $$\theta(z) \psi(z) \theta(z) = \theta(z), \quad \mbox{for}\; z \in \mathbb{B}^m.$$
\end{Corollary}

In conclusion, Corollary \ref{main-cor} shows that Statement \ref{Q4} implies Statement \ref{Q3}.

\newsection{Resolutions of Hilbert modules over $\mathbb{C}[z_1, \ldots, z_m]$}

Consideration of resolutions such as those in the preceding section and the ones given in Theorem \ref{AMVMT} raises the question of what kind of resolutions exist for pure co-spherically contractive Hilbert modules over $\mathbb{C}[z_1, \ldots, z_m]$. In particular, Theorem \ref{AMVMT} yields a unique resolution of an arbitrary pure co-spherically contractive Hilbert module $\clm$ over $\mathbb{C}[z_1, \ldots, z_m]$ in terms of DA-modules and inner multipliers. More specifically, consider DA-modules $\{H^2_m \otimes \cle_k\}$ for Hilbert spaces $\{\cle_k\}$ and inner multipliers $\varphi_k \in \clm(\cle_k, \cle_{k-1})$, or partially isometric module maps $\{M_{\varphi_k}\}$ for $k \geq 1$; set $$X_k = M_{\varphi_k} : H^2_m \otimes \cle_k \raro H^2_m \otimes \cle_{k-1}, \;\; k \geq 1;$$ and a co-isometric module map $$X_0 = \pi_{\clm} : H^2_m \otimes \cle_0 \raro \clm,$$ which is exact. That is, $\mbox{ran~} X_{k+1} = \mbox{ker~} X_k$ for $k \geq 1$. Here $k = 0, 1, \ldots, N$, with the possibility of $N = + \infty$. A basic question is whether such a resolution can have finite length or, equivalently, whether we can take $\cle_N = \{0\}$ for some finite $N$. That will be the case if and only if some $X_k$ is an isometry or, equivalently, if $\mbox{ker} X_k = \{0\}$. Unfortunately, the following result shows that this is not possible when $m >1$, unless $\clm$ is a DA-module and the resolution is a trivial one.

\begin{Theorem}\label{res}
For $m>1$, if $V : H^2_m \otimes \cle \raro H^2_m \otimes \cle_*$ is an isometric module map for Hilbert spaces $\cle$ and $\cle_*$, then there exists an isometry $V_0 : \cle \raro \cle_*$ such that $$V( \bm{z}^{\bm{k}} \otimes x)  = \bm{z}^{\bm{k}} \otimes V_0 x, \;\; \mbox{for~} \bm{k} \in \mathbb{N}^m, x \in \cle_*.$$Moreover, $\mbox{ran~} V$ is a reducing submodule of $H^2_m \otimes \cle_*$ of the form $H^2_m \otimes (\mbox{ran}\, V_0)$.
\end{Theorem}
\NI \textsf{Proof.} For $x \in \cle$, $\|x\| = 1$, we have $$V (1 \otimes x) = f(\bm{z}) = \sum_{\bm{k} \in \mathbb{N}^m} a_{\bm{k}} \bm{z}^{\bm{k}}, \quad \quad \mbox{for} \, \{a_k\} \subseteq \cle.$$ Then $$V (z_1 \otimes x) = V M_{z_1} (1 \otimes x) = M_{z_1} V (1 \otimes x) = M_{z_1} f = z_1 f,$$ and $$\|z_1 f \|^2 = \|z_1 V(1 \otimes x)\|^2 = \|z_1 \otimes x\|^2 = 1 = \|f\|^2.$$Therefore, we have $$\sum_{\bm{k} \in \mathbb{N}^m} \|a_{\bm{k}}\|^2_{\cle_*} \|\bm{z}^{\bm{k}}\|^2 = \sum_{\bm{k} \in \mathbb{N}^m} \|a_{\bm{k}}\|^2_{\cle_*} \|\bm{z}^{\bm{k}+e_1}\|^2, \; \; \mbox{where~} \bm{k} + e_1 = (k_1+1, \ldots, k_m),$$ or $$\sum_{\bm{k} \in \mathbb{N}^m} \|a_{\bm{k}}\|^2_{\cle_*} \{ \|\bm{z}^{\bm{k} + e_1}\|^2 - \|\bm{z}^{\bm{k}}\|^2\} = 0.$$
If $\bm{k} = (k_1, \ldots k_m)$, then
\[
\begin{split}
\|\bm{z}^{\bm{k} + e_1}\|^2 & = \frac{(k_1 + 1)! \cdots k_m!}{(k_1 + \cdots + k_m +1)!} = {\frac{k_1! \cdots k_m!}{(k_1 + \cdots + k_m)!}} \frac{k_1 + 1}{k_1 + \cdots + k_m + 1}\\ & < {\frac{k_1! \cdots k_m!}{(k_1 + \cdots + k_m)!}} = \|\bm{z}^{\bm{k}}\|^2,
\end{split}
\]
unless $k_2 = k_3 = \ldots = k_m = 0$. Since, $a_{\bm{k}} \neq 0$ implies $\|\bm{z}^{\bm{k}+e_1}\| = \|\bm{z}^{\bm{k}}\|$ we have  $k_2 = \cdots = k_m = 0$. Repeating this argument using $i = 2, \ldots,m$, we see that $a_{\bm{k}} = 0$ unless $\bm{k} = (0, \ldots, 0)$ and therefore, $f(\bm{z}) = 1 \otimes y$ for some $y \in \cle_*$. Set $V_0 x = y$ to complete the first part of the proof.

Finally, since $\mbox{ran~} V = H^2_m \otimes (\mbox{ran} V_0)$, we see that $\mbox{ran~} V$ is a reducing submodule, which completes the proof. \qed

Note that this result generalizes Corollary 3.3 of \cite{DS} and is related to an earlier result of Guo, Hu and Xu \cite{GHX}.

The theorem implies that all resolutions by DA-modules with partially isometric maps are trivial in a sense we will make precise. We start with a definition.

\begin{Definition}
An inner resolution of length $N$, for $N = 1, 2, 3, \ldots, \infty$, for a pure co-spherical contractive Hilbert module $\clm$ is given by a collection of Hilbert spaces $\{\cle_k\}_{k=0}^N$, inner multipliers $\varphi_k \in \clm(\cle_k, \cle_{k-1})$ for $k = 1, \ldots, N$ with $X_k = M_{\varphi_k}$ and a co-isometric module map $X_0 : H^2_m \otimes \cle_0 \raro \clm$ so that $$\mbox{ran} \, X_k = \mbox{ker}\, X_{k-1},$$ for $k = 0, 1, \ldots, N$. To be more precise, for $N < \infty$ one has the finite resolution $$ 0 \longrightarrow H^2_m \otimes \cle_N \stackrel{X_N}{\longrightarrow} H^2_m \otimes \cle_{N-1} \longrightarrow \cdots \longrightarrow H^2_m \otimes \cle_1 \stackrel{X_1} \longrightarrow H^2_m \otimes \cle_0 \stackrel{X_0} \longrightarrow \clm \longrightarrow 0,$$ and for $N = \infty$, the infinite resolution
$$\cdots \longrightarrow H^2_m \otimes \cle_N \stackrel{X_{N}}{\longrightarrow} H^2_m \otimes \cle_{N-1} \longrightarrow \cdots \longrightarrow H^2_m \otimes \cle_1 \stackrel{X_{1}} \longrightarrow H^2_m \otimes \cle_0 \stackrel{X_{0}} \longrightarrow \clm \longrightarrow 0.$$
\end{Definition}

\begin{Theorem}\label{trivial-res}
If the pure, co-spherically contractive Hilbert module $\clm$ possesses a finite inner resolution, then $\clm$ is isometrically isomorphic to $H^2_m \otimes \clf$ for some Hilbert space $\clf$.
\end{Theorem}

\NI\textsf{Proof.} Applying the previous theorem to $M_{\varphi_N}$, we decompose $\cle_{N-1} = \cle^1_{N-1} \oplus \cle^2_{N-1}$ so that $\tilde{M}_{\psi_{N-1}} = M_{N-1}|_{H^2_m \otimes \cle^2_{N-1}} \in \cll(H^2_m \otimes \cle^2_{N-1}, H^2_m \otimes \cle_{N-2})$ is an isometry onto $\mbox{ran} \, M_{N-1}$. Hence, we can apply the theorem to $\tilde{M}_{N-1}$. Therefore, using induction we obtain the desired conclusion. \qed

The following statement proceeds directly from the theorem.

\begin{Corollary}
If $\theta\in \clm(\cle, \cle_*)$ is an inner multiplier for the Hilbert spaces $\cle$ and $\cle_*$ with $\mbox{ker~}M_{\theta} = \{0\}$, then the quotient module $\clh_{\theta} = (H^2_m \otimes \cle_*)/ \, \mbox{ran} \,M_{\theta}$ is isometrically isomorphic to $H^2_m \otimes \clf$ for  a Hilbert space $\clf$. Moreover, $\clf$ can be identified with $(\mbox{ran}\, V_0)^{\perp}$, where $V_0$ is the isometry from $\cle$ to $\cle_*$ given in Theorem \ref{res}.
\end{Corollary}

Note that in the preceding corollary, one has $\mbox{dim~} \cle_* = \mbox{dim~} \cle + \mbox{dim~}\clf$.

A resolution of $\clm$ can always be made longer in a trivial way. Suppose we have the resolution $$ 0 \longrightarrow H^2_m \otimes \cle_N \stackrel{X_N}{\longrightarrow} H^2_m \otimes \cle_{N-1} \longrightarrow \cdots \longrightarrow H^2_m \otimes \cle_0 \stackrel{X_0}\longrightarrow \clm \longrightarrow 0.$$ If $\cle_{N+1}$ is a nontrivial Hilbert space, then define $X_{N+1}$ as the inclusion map of $H^2_m \otimes \cle_{N+1} \subseteq H^2_m \otimes (\cle_N \oplus \cle_{N+1})$. Further, set $\tilde{X}_N$ equal to $X_N$ on $H^2_m \otimes \cle_{N} \subseteq H^2_m \otimes (\cle_{N+1} \oplus \cle_N)$ and equal to $0$ on $H^2_m \otimes \cle_{N+1} \subseteq H^2_m \otimes (\cle_N \oplus \cle_{N+1})$. Extending $\tilde{X}_{N}$ to all of $H^2_m \otimes \cle_{N+1}$ linearly, we obtain a longer resolution essentially equivalent to the original one $$ 0 \longrightarrow H^2_m \otimes \cle_{N+1} \stackrel{X_{N+1}}{\longrightarrow} H^2_m \otimes (\cle_{N+1} \oplus \cle_N)  \stackrel{\tilde{X}_N} \longrightarrow \cdots \longrightarrow  \clm \longrightarrow 0.$$
Moreover, the new resolution will be inner if the original one is.

The proof of the preceding theorem shows that any finite inner resolution by DA-modules is equivalent to a series of such trivial extensions of the resolution $$ 0 \longrightarrow H^2_m \otimes \cle \stackrel{X}{\longrightarrow} H^2_m \otimes \cle \longrightarrow 0,$$ for some Hilbert space $\cle$ and $X = I_{H^2_m \otimes \cle}$. We will refer to such a resolution as a {\it trivial inner resolution}. We use that terminology to summarize this supplement to the theorem in the following statement.

\begin{Corollary}
All finite inner resolutions for a pure co-spherically contractive Hilbert module $\clm$ are trivial inner resolutions.
\end{Corollary}

What happens when we relax the conditions on the module maps $\{X_k\}$ so that $\mbox{ran}\, X_{k} = \mbox{ker}\, X_{{k-1}}$ for all $k$ but do not require them to be partial isometries? In this case, non-trivial finite resolutions do exist, completely analogous to what happens for the case of the Hardy or Bergman modules over $\mathbb{C}[z_1, \ldots, z_m]$ for $m>1$. We describe a simple example.

Consider the module $\mathbb{C}_{(0,0)}$ over $\mathbb{C}[z_1, z_2]$ defined so that $$p(z_1, z_2) \cdot \lambda = p(0,0) \lambda, \;\;\mbox{where~} p \in \mathbb{C}[z_1, z_2] \,\mbox{and~} \lambda \in \mathbb{C},$$ and the resolution:

$$ 0 \longrightarrow H^2_2 \stackrel{X_2}{\longrightarrow} H^2_2 \oplus H^2_2 \stackrel{X_1} \longrightarrow H^2_2 \stackrel{X_0} \longrightarrow \mathbb{C}_{(0,0)} \longrightarrow 0,$$
where $X_0 f = f(0,0)$ for $f \in H^2_2$, $X_1(f_1 \oplus f_2) = M_{z_1} f_1 + M_{z_2} f_2$ for $f_1 \oplus f_2 \in H^2_2 \oplus H^2_2$, and $X_2 f = M_{z_2} f \oplus ( - M_{z_1} f)$ for $f \in H^2_2$. One can show that this sequence, which is closely related to the Koszul complex, is exact and non-trivial; in particular, it does not split as trivial resolutions do.

Another question one can ask is the relationship between the inner resolution for a pure co-spherically contractive Hilbert module and more general, {\it not necessarily inner}, resolutions by DA-modules. In particular, is there any relation between the minimal length of a not necessarily inner resolution and the inner resolution. Theorem \ref{incl} provides some information on this matter.

A parallel notion of resolution for Hilbert modules was studied by Arveson \cite{A04}, which is different from the one considered in this paper. For Arveson, the key issue is the behavior of the resolution at $\bm{0} \in \mathbb{B}^m$ or the localization of the sequence of connecting maps at $\bm{0}$. His main goal, which he accomplishes and is quite non trivial, is to extend an analogue of the Hilbert's syzygy theorem. In particular, he exhibits a resolution of Hilbert modules in his class which ends in finitely many steps. The resolutions considered in (\cite{DM1}, \cite{DM2}) and this paper are related to dilation theory although the requirement that the connecting maps are partial isometries is sometimes relaxed.

\newsection{Hilbert modules over $\mathbb{F}[Z_1, \ldots, Z_m]$}

Although we use the following lemma only in the non-commutative case, it also holds in the commutative case as indicated.

\begin{Lemma}\label{LM1}
If $\clh$ is a co-spherically contractive Hilbert module over $\mathbb{C}[z_1, \ldots, z_m]$ or $\mathbb{F}[Z_1, \ldots, Z_m]$, respectively, which is similar to $H^2_m \otimes \clf$, or $F^2_m \otimes \clf$, respectively, for some Hilbert space $\clf$, then $\clh$ is pure.
\end{Lemma}
\NI \textsf{Proof.} We use the notation for the commutative case but the proof in both cases is the same. Let $X : \clh \rightarrow H^2_m \otimes \clf$ be an invertible module map. Then $M_i = X^{-1} M_{z_i} X$ for all $i = 1, \ldots, m$. Since $\{P^l_{\clh} (I_{\clh})\}^{\infty}_{l=0}$ is a decreasing sequence of positive operators, it suffices to show that $$\mbox{WOT} - \mbox{lim}_{l \raro \infty} P_{\clh}^l (I_{\clh}) = 0.$$ To see that this is the case, let $f_1$ be a vector in $\clh$ and set $f = X^{*-1} f_1$. Then
\begin{equation*}
\begin{split}
 \la \sum_{|\bm{k}| = l} M^{\bm{k}} M^{*\bm{k}} f_1, f_1 \ra  = & \la \sum_{|\bm{k}| = l} X^{-1} M_z^{\bm{k}} X X^* M_z^{*\bm{k}} X^{*-1} f_1, f_1 \ra =  \la\sum_{|\bm{k}|=l} M_z^{\bm{k}} X X^* M_z^{*\bm{k}} f, f\ra \\ & \leq \|X\|^2 \sum_{|\bm{k}|=l} \langle M_z^k M^{*k}_z f, f\rangle.
\end{split}
\end{equation*}
\NI Letting $l \raro \infty$ in the last expression, we conclude that the required limit is zero, which completes the proof.\qed

\vspace{0.1in}

Actually, the proof shows that two similar co-spherically contractive Hilbert modules over $\mathbb{C}[z_1, \ldots, z_m]$, or two similar contractive Hilbert modules over $\mathbb{F}[Z_1, \ldots, Z_m]$, are either both pure or both not pure.

\begin{Theorem}\label{TH2}
Let $\clh$ be a pure co-spherically contractive Hilbert module over $\mathbb{F}[Z_1, \ldots, Z_m]$. Then $\clh$ is similar to $F^2_m \otimes \clf$ for some Hilbert space $\clf$ if and only if the characteristic operator $\Theta$ of $\clh$ in $\cll(F^2_m \otimes \cle, F^2_m \otimes \cle_*)$, for some Hilbert spaces $\cle$ and $\cle_*$, is left invertible; that is, if and only if there exists a multi-analytic operator $\Psi : F^2_m \otimes {\cle_*} \raro F^2_m \otimes \cle$ such that $$\Psi \Theta = I_{F^2_m \otimes \cle}.$$
\end{Theorem}

\NI \textsf{Proof.} First, using Theorem \ref{Pop1} we realize the pure contractive Hilbert module $\clh$ as the quotient module given by its characteristic function $\Theta$, which is an isometric multi-analytic map. That is,  $$\clh \cong \clh_{\Theta} = (F^2_m \otimes \cle_*)/\Theta(F^2_m \otimes \cle).$$ Now given a module map $X : F^2_m \otimes \clf \raro \clh_{\Theta} $, we appeal to the noncommutative analogue of the CLT (see Theorem 6.1 in \cite{P89} or Theorem 5.1 in \cite{P06}) to obtain a multi-analytic operator $\Phi : F^2_m \otimes \clf \raro F^2_m \otimes \cle_*$ such that $$P_{\clh_{\Theta}} \Phi = X.$$ Consider, the bounded module map $$Z :(F^2_m \otimes \clf) \oplus (F^2_m \otimes \cle) \rightarrow F^2_m \otimes \cle_*$$ defined by $$Z(f \oplus g) = \Phi f + \Theta g,$$ for all $f \oplus g \in (F^2_m \otimes \clf) \oplus (F^2_m \otimes \cle)$. Then $Z$ is invertible if and only if $X$ is invertible. This follows by noting that $X$ is invertible if and only if the range of $Z$, which is the span of $\clh_{\Theta}$ and $\mbox{ran}\, \Theta$ is $F^2_m \otimes \cle_*$, and $X$ is one-to-one if and only if $Z$ is.

To prove the necessity part of the theorem, assume that $X$ is invertible or, equivalently, that $Z$ is invertible. Consequently, we can define a module idempotent $Q$ on $F^2_m \otimes \cle_*$ such that $$Q \Theta = \Theta$$ and $$\mbox{ran}\, Q = \mbox{ran}\, \Theta.$$ Then the bounded module map $\hat{Q} : F^2_m \otimes \cle_* \rightarrow F^2_m \otimes \cle$ defined by $$\hat{Q} (\Phi f + \Theta g) = g, \quad \quad \quad \quad \quad \Phi f + \Theta g \in (F^2_m \otimes \cle_*)$$ satisfies $$Q = \Theta \hat{Q}.$$
Since $\hat{Q}$ is a module map, there exists a multi-analytic operator $\Psi : F^2_m \otimes \cle_* \rightarrow F^2_m \otimes \cle$ such that $$\hat{Q} = \Psi.$$ Hence $$\Theta = Q \Theta = \Theta \hat{Q} \Theta = \Theta \Psi \Theta.$$ Since $\Theta$ is an isometry, the necessity  part follows; that is, $\Theta$ has a left inverse.

To prove the sufficiency part, let $\Psi : F^2_m \otimes \cle_* \raro F^2_m \otimes \cle$ be a multi-analytic operator such that $$\Psi \Theta = I_{F^2_m \otimes \cle}.$$ Then $Q = \Theta \Psi$ is an idempotent on $F^2_m \otimes \cle_*$ and any $f$ in $F^2_m \otimes \cle_*$ can be expressed as $$f = (f - \Theta \Psi f) + \Theta \Psi f,$$ where $f - \Theta \Psi f$ is in $\mbox{ker} \, \Psi$ and $\Theta \Psi f$ is in $\mbox{ran} \, \Theta$. Thus, $$\mbox{ran}\, Q = \mbox{ran}\, \Theta, \quad \quad \mbox{and} \quad \quad \mbox{ker}\, \Psi = \mbox{ran} \, (I - Q).$$ Since $\mbox{ker}\, \Psi$ is a submodule of $F^2_m \otimes \cle_*$, by Theorem \ref{BLH}, the noncommutative version of the BLHT, there exists an inner multi-analytic operator $\Phi : F^2_m \otimes \clf \rightarrow F^2_m \otimes \cle_*$ for some Hilbert space $\clf$ such that $$\mbox{ker~} \Psi = \mbox{ran~} (I - Q) = \Phi (F^2_m \otimes \clf),$$  Consequently, $$F^2_m \otimes \cle_* = \mbox{ran~} \Phi \stackrel{\bm{.}}{+} \mbox{ran~} \Theta.$$ Then one can define the invertible module map $Z$ as in the necessity part. Setting $X = P_{\clh_{\Theta}} \Phi$ defines the required similarity between $\clh_{\Theta}$ and $F^2_m \otimes \clf$, which completes the proof. \qed

As mentioned in the introduction, specializing the preceding proof to the (commutative) $m=1$ case yields a new proof of the old result on the similarity of contraction operators to unilateral shifts.

The main difference in the above proof and that of Corollary \ref{main-cor} for the commutative case is that here we can assume that $\Theta$ has no kernel and one of the complemented submodule is isomorphic to a DA-module.

In the proof of Theorem \ref{TH2}, we did not use the fact that the characteristic function is an isometry but the fact that $\mbox{ker~}\Theta = \{0\}$ and $\mbox{ran~}\Theta$ is closed. Hence we can state a more general result in terms of a module resolution.

\begin{Theorem}
Let $\cle$ and $\cle_*$ be Hilbert spaces and $\Theta : F^2_m \otimes \cle \raro F^2_m \otimes \cle_*$ be a multi-analytic operator such that $\mbox{ker~}\Theta = \{0\}$ and $\mbox{ran~}\Theta$ is closed. Then the quotient space $\clh_{\Theta}$, given by $(F^2_m \otimes \cle_*)/\, \mbox{ran}\, \Theta$ is similar to $F^2_m \otimes \clf$ for some Hilbert space $\clf$ if and only if $\Theta \Psi \Theta = \Theta$ for some multi-analytic operator $\Psi : F^2_m \otimes \cle_* \raro F^2_m \otimes \cle$.
\end{Theorem}

\vspace{0.2in}

\newsection{Concluding remarks}

Observe that if $\clh$ is a Hilbert module over $\mathbb{C}[z_1, \ldots, z_m]$ (or $A(\Omega)$, where $\Omega$ is a bounded connected open subset of $\mathbb{C}^m$), Corollary \ref{main-cor} remains true under the assumption that the analogue of the CLT holds for the class of Hilbert modules under consideration. In particular, Corollary \ref{main-cor} can be generalized to any reproducing kernel Hilbert module where the kernel is given by a complete Nevanlinna-Pick kernel.

\end{document}